\documentclass[11pt]{article}

\usepackage{scrpage2}
\clearscrheadfoot
\usepackage{tikz}
\usetikzlibrary{angles, quotes, babel}
\usetikzlibrary{decorations.pathreplacing}
\usetikzlibrary{spy}
\tikzset{SpyStyle/.style={spy using outlines={rectangle, magnification=15, width=2cm, height=2cm, connect spies, blue!70!black}}}

\usepackage{pgfplots}
\pgfplotsset{compat=1.13}
\usepgfplotslibrary{patchplots}
\definecolor{Farbe}{rgb}{0.8, 0.8, 0.8}

\usepackage{caption}
\usepackage{subcaption}

\usepackage{times}
\usepackage{amsmath}
\usepackage{amsfonts}
\usepackage{amssymb}

\usepackage{tabularx}
\usepackage[linkcolor=black,urlcolor=black,citecolor=black]{hyperref}
\usepackage[all]{hypcap}
\hypersetup{linktocpage,colorlinks=true,pdfborder={0 0 0}}

\usepackage{geometry}
\title{\textsc{A 3-regular matchstick graph of girth 5 consisting of 54 vertices}}

\author{
	\\ \Large{Mike Winkler$^\star$\quad Peter Dinkelacker$^{\star\star}$\quad Stefan Vogel$^{\star\star\star}$}
}

\date{
	\small$^\star$Fakult\"at f\"ur Mathematik, Ruhr-Universit\"at Bochum, Germany, mike.winkler@ruhr-uni-bochum.de\\[2mm]
	$^{\star\star}$Togostr. 79, 13351 Berlin, Germany, peter@grity.de\\[2mm]
	$^{\star\star\star}$Raun, Dorfstr. 7, 08648 Bad Brambach, Germany, backebackekuchen16@gmail.com\\[10mm]
	\large{November 5, 2019}
}

\begin{document}
  
  \maketitle
  
  \begin{abstract}
    In 2010 it was proved that a 3-regular matchstick graph of girth 5 must consist at least of 30 vertices. The smallest known example consisted of 180 vertices. In this article we construct an example consisting of 54 vertices and prove its geometrical correctness.\\[10mm]
  \end{abstract}
  
  \begin{center}
  	\Large\textsc{Introduction}
  \end{center}
  
  A 3-regular matchstick graph is a planar unit-distance graph whose vertices have all degree 3. The girth of a graph is the length of a shortest cycle contained in the graph. Therefore, no rigid matchstick graph of girth $\geq4$ exists, because such graphs contain only flexible subgraphs.
  
  In 2010 Sascha Kurz and Giuseppe Mazzuoccolo proved that a 3-regular matchstick graph of girth 5 consists at least of 30 vertices and gave an example consisting of 180 vertices \cite{Kurz}.
  
  In this paper we construct an example with 54 vertices (see Figure 1) and prove its geometrical correctness. This graph was first presented by the authors on February 8, 2019 in a graph theory internet forum \cite{Thread}. The graph and its proof by construction can also be viewed on the authors website \href{http://mikematics.de/matchstick-graphs-calculator.htm}{\textit{mikematics.de}}\footnote{http://mikematics.de/matchstick-graphs-calculator.htm} with a software called \textsc{Matchstick Graphs Calculator} (MGC) \cite{MGC}. This software runs directly in web browsers.\footnote{For optimal functionality and design please use the Firefox web browser.} The MGC includes an animation function for representing the movement we are mentioning in our proof.
  
  We also found further 3-regular matchstick graphs of girth 5 with less than 70 vertices consisting of 58, 60, 64, 66, and 68 vertices. These graphs are included in the MGC and exhibit in a separate paper \cite{Winkler}.
  
  \newpage
  
  \begin{center}
  	\Large\textsc{The construction of the graph}
  \end{center}
  \quad
  
  \textbf{Theorem.} \textit{A 3-regular matchstick graph of girth 5 consisting of 54 vertices exists.}
  \\
  
  \textit{Proof.} Given are the vertices $P1$ and $P2$ with the distance of a unit length in the plane.
  In the following an 'edge' means always an 'edge of unit length'.
  Add an edge from $P2$ to $P1$.
  Add an edge from $P3$ to $P1$ with angle $\alpha=102^\circ$ to the edge from $P2$ to $P1$.
  Add an edge from $P4$ to $P3$ with angle $\beta=67^\circ$ to the edge from $P1$ to $P3$.
  Add an edge from $P5$ to $P3$ with angle $\gamma=74^\circ$ to the edge from $P4$ to $P3$.
  Add an edge from $P6$ to $P5$ with angle $\delta=81^\circ$ to the edge from $P3$ to $P5$.
  Add two edges to complete the isosceles triangle $P7,P2,P4$, the base remains open.
  Add two edges to complete the isosceles triangle $P8,P4,P6$, the base remains open.
  Add an edge from $P9$ to $P8$ with angle $\epsilon=24^\circ$ to the edge from $P4$ to $P8$.
  Add two edges to complete the isosceles triangle $P10,P9,P7$, the base remains open.
  Add an edge from $P11$ to $P5$ with angle $\zeta=69^\circ$ to the edge from $P6$ to $P5$.
  Add an edge from $P12$ to $P11$ with angle $\eta=106^\circ$ to the edge from $P5$ to $P11$.
  Add two edges to complete the isosceles triangle $P13,P6,P12$, the base remains open.
  Add two edges to complete the isosceles triangle $P14,P13,P9$, the base remains open.
  Add an edge from $P15$ to $P14$ with angle $\theta=3^\circ$ to the edge from $P9$ to $P14$.
  Add an edge from $P16$ to $P11$ with angle $\iota=61^\circ$ to the edge from $P12$ to $P11$.
  Add two edges to complete the isosceles triangle $P17,P12,P15$, the base remains open.
  Add two edges to complete the isosceles triangle $P18,P17,P16$, the base remains open.
  Add an edge from $P19$ to $P1$ with angle $\kappa=85^\circ$ to the edge from $P2$ to $P1$.
  Add an edge from $P20$ to $P19$ with angle $\lambda=91^\circ$ to the edge from $P1$ to $P19$.
  Add two edges to complete the isosceles triangle $P21,P20,P2$, the base remains open.
  Add two edges to complete the isosceles triangle $P22,P10,P21$, the base remains open.
  Add two edges to complete the isosceles triangle $P23,P15,P22$, the base remains open.
  Add two edges to complete the isosceles triangle $P24,P18,P23$, the base remains open.
  Add an edge from $P25$ to $P19$ with angle $\mu=38^\circ$ to the edge from $P20$ to $P19$.
  Add an edge from $P26$ to $P16$ with angle $\nu=65^\circ$ to the edge from $P18$ to $P16$.
  Add two edges to complete the isosceles triangle $P27,P24,P26$, the base remains open.
  Connect $P25$ and $P26$ by a copy of the subgraph $P1-P27$. The copied subgraph matches the already existing graph, because of congruence.
  Add two edges to complete the isosceles triangle $P53,P47,P27$, the base remains open.
  Add two edges to complete the isosceles triangle $P54,P20,P52$, the base remains open.
  Add an edge of unknown length from $P53$ to $P54$.
  
  We get a point-symmetric 3-regular planar graph of girth 5 consisting of 54 vertices. The graph contains clearly visible distances between vertices and edges (see Figure 1). The method of construction ensures that each edge, except $P53,P54$, has unit length. The distance $P53,P54$ measures $\approx1.0007$. Varying the angle $\mu$ has no effect on the other twelve angles mentioned before. Because of the point symmetry we can separate the graph into two rigid subgraphs $\mathcal{G}_1$ and $\mathcal{G}_2$. $\mathcal{G}_1$ consisting of $P1-P24,P26,P27$, and $\mathcal{G}_2$ consisting of $P25,P28-P52$. By holding $\mathcal{G}_1$ fixed and changing $\mu$ continuously close to its value from about $39^\circ$ to $37^\circ$, the distance $P53,P54$ varies around the unit length from approximately $0{,}991$ to $1{,}012$. During this changing $\mathcal{G}_2$ moves slightly from the right to the left, while $P53$ moves up and $P54$ moves down. During the movement there are no unauthorized overlaps or contacts. This shows that a 3-regular matchstick graph of girth 5 consisting of 54 vertices exists. The distance $P53,P54$ has unit length, if $\mu\approx38{,}067338069376^\circ$.\hfill$\square$
  
  \newpage
  
  \begin{center}
  	\begin{minipage}{\linewidth}
  		\centering
  		
  		\xdef\LstPN{0}
  		\newif\ifDupe
  		\pgfplotsset{avoid dupes/.code={\Dupefalse
  				\xdef\anker{\DefaultTextposition} 
  				\foreach \X in \LstPN
  				{\pgfmathtruncatemacro{\itest}{ifthenelse(\X==\punktnummer,1,0)}
  					\ifnum\itest=1
  					\global\Dupetrue
  					\breakforeach
  					\fi}
  				\ifDupe
  				\typeout{\punktnummer\space ist\space ein\space Duplikat!}%
  				\xdef\punktnummer{} 
  				\else
  				\xdef\LstPN{\LstPN,\punktnummer}
  				\typeout{\punktnummer\space ist\space neu\space mit\space urprgl.\space Anker=\anker}
  				\foreach \X in \LstExcept
  				{\ifnum\X=\punktnummer
  					\xdef\anker{\AusnahmeTextposition}
  					\fi}
  				\typeout{\punktnummer\space ist\space neu\space mit\space Anker=\anker}
  				\fi}}
  		
  		\def\DefaultTextposition{east}
  		\def\AusnahmeTextposition{west}
  		\def\AusnahmeListe{}
  		\xdef\BeliebigesVorhandenesKoordinatenpaar{{3.37452506686058084640, 0.62932039104983827915}} 
  		\colorlet{Kantenfarbe}{gray}
  		\colorlet{Punktfarbe}{red}
  		\def\Beschriftung{} 
  		\pgfplotsset{
  			x=30.0mm, y=30.0mm,  
  		}
  		
  		\xdef\LstExcept{\AusnahmeListe}
  		\pgfdeclarelayer{bg}    
  		\pgfsetlayers{bg,main}  
  		
  		\pgfmathsetmacro{\xAlias}{\BeliebigesVorhandenesKoordinatenpaar[0]}
  		\pgfmathsetmacro{\yAlias}{\BeliebigesVorhandenesKoordinatenpaar[1]}
  		
  		\begin{tikzpicture}[SpyStyle]
  		
  		\begin{axis}[hide axis, colormap={kantenfarbe}{color=(Kantenfarbe) color=(Kantenfarbe)}, line width=0.2, 
  		]
  		\addplot+[mark size=0.8pt, 
  		mark options={Punktfarbe}, 
  		table/row sep=newline, 
  		patch, 
  		patch type=polygon,
  		vertex count=2, 
  		%
  		patch table with point meta={
  			Startpkt Endpkt colordata  \\
  			1 1 \\
  			2 1 \\
  			3 1 \\
  			4 3 \\
  			5 3 \\
  			6 5 \\
  			7 2 \\
  			7 4 \\
  			8 4 \\
  			8 6 \\
  			9 8 \\
  			10 9 \\
  			10 7 \\
  			11 5 \\
  			12 11 \\
  			13 6 \\
  			13 12 \\
  			14 13 \\
  			14 9 \\
  			15 14 \\
  			16 11 \\
  			17 12 \\
  			17 15 \\
  			18 17 \\
  			18 16 \\
  			19 1 \\
  			20 19 \\
  			21 20 \\
  			21 2 \\
  			22 10 \\
  			22 21 \\
  			23 15 \\
  			23 22 \\
  			24 18 \\
  			24 23 \\
  			25 19 \\
  			25 43 \\
  			26 16 \\
  			26 46 \\
  			27 24 \\
  			27 26 \\
  			28 28 \\
  			29 28 \\
  			30 28 \\
  			31 30 \\
  			32 30 \\
  			33 32 \\
  			34 29 \\
  			34 31 \\
  			35 31 \\
  			35 33 \\
  			36 35 \\
  			37 34 \\
  			37 36 \\
  			38 32 \\
  			39 38 \\
  			40 33 \\
  			40 39 \\
  			41 36 \\
  			41 40 \\
  			42 41 \\
  			43 38 \\
  			44 39 \\
  			44 42 \\
  			45 43 \\
  			45 44 \\
  			46 28 \\
  			47 46 \\
  			48 29 \\
  			48 47 \\
  			49 37 \\
  			49 48 \\
  			50 42 \\
  			50 49 \\
  			51 45 \\
  			51 50 \\
  			52 51 \\
  			52 25 \\
  			53 47 \\
  			53 27 \\
  			53 54 \\
  			54 20 \\
  			54 52 \\
  		},
  		%
  		visualization depends on={value \thisrowno{0} \as \punktnummer},
  		every node near coord/.append style={
  			/pgfplots/avoid dupes,
  		},
  		nodes near coords={\Beschriftung},
  		nodes near coords style={
  			anchor=\anker,
  			text=black, font=\scriptsize, 
  			name=p-\punktnummer, 
  			path picture={
  				\coordinate[] (P\punktnummer) at (p-\punktnummer.\anker);}
  		},
  		]
  		table[header=true, x index=1, y index=2, row sep=\\] {
  			Nr x y        \\
  			0 0 0         \\
  			1 3.37452506686058084640 0.62932039104983827915  \\
  			2 2.92053456712103409743 1.52032691523820595592  \\
  			3 2.59737910540361038869 0.00000000000000074670  \\
  			4 2.32174174958661128088 0.96126169593831956028  \\
  			5 1.59737910540361016665 0.00000000000000000000  \\
  			6 1.75381357044384067478 0.98768834059513788137  \\
  			7 3.24370588969779527488 0.57398643937130799930  \\
  			8 1.99321392857307722046 0.01676740224454608200  \\
  			9 1.90917858874795776103 1.01323017709390117425  \\
  			10 2.79894932374493476956 1.46963793244314944175  \\
  			11 0.73135370161917112597 0.49999999999999916733  \\
  			12 0.97327559721883805111 1.47029572627599591783  \\
  			13 0.89627750840113462605 0.47326448588358921432  \\
  			14 1.01749550265636368174 1.46589039621876637831  \\
  			15 1.93164697348957514755 1.06051761916846154499  \\
  			16 0.00000000000000000000 1.18199836006249725529  \\
  			17 1.78799838453970760810 2.05014637930668230936  \\
  			18 0.94252918542962682569 1.51612219131876080347  \\
  			19 4.22257316301700669214 1.15923965528304373507  \\
  			20 3.70753508810695242559 2.01640695598515540254  \\
  			21 3.78607992642551360873 1.01949637405039128524  \\
  			22 3.64107028624210427026 2.00892661628774993332  \\
  			23 2.68390421303229942396 1.71938714005674819774  \\
  			24 1.75742411389290520418 2.09573090033272002586  \\
  			25 4.34560893271744674848 2.15164189196644661806  \\
  			26 0.09551101588407498455 2.17742673309834344764  \\
  			27 0.89922600971951316406 1.58241227219499536183  \\
  			28 1.06659488174094096991 3.69974823401495145347  \\
  			29 1.52058538148048727479 2.80874170982658366569  \\
  			30 1.84374084319791275988 4.32906862506479050978  \\
  			31 2.11937819901491186769 3.36780692912647072745  \\
  			32 2.84374084319791231579 4.32906862506479050978  \\
  			33 2.68730637815768247378 3.34138028446965229534  \\
  			34 1.19741405890372765164 3.75508218569348128923  \\
  			35 2.44790602002844526197 4.31230122282024286307  \\
  			36 2.53194135985356494345 3.31583844797088955758  \\
  			37 1.64217062485658726878 2.85943069262164062394  \\
  			38 3.70976624698235113442 3.82906862506479006569  \\
  			39 3.46784435138268465337 2.85877289878879459195  \\
  			40 3.54484244020038818945 3.85580413918120079586  \\
  			41 3.42362444594515880070 2.86317822884602346534  \\
  			42 2.50947297511194733488 3.26855100589632918684  \\
  			43 4.44111994860152226039 3.14707026500229281041  \\
  			44 2.65312156406181465229 2.27892224575810820042  \\
  			45 3.49859076317189554572 2.81294643374602948427  \\
  			46 0.21854678558451454129 3.16982896978174766289  \\
  			47 0.73358486049456927969 2.31266166907963333088  \\
  			48 0.65504002217600776348 3.30957225101440011272  \\
  			49 0.80004966235941687991 2.32014200877703968828  \\
  			50 1.75721573556922283643 2.60968148500804275614  \\
  			51 2.68369583470861705621 2.23333772473206959575  \\
  			52 3.54189393888200898530 2.74665635286979448182  \\
  			53 1.72070107156966001227 2.15265671342901310226  \\
  			54 2.72041887703186180403 2.17641191163577696344  \\
  		};
  		
  		\addplot[no marks, 
  		nodes near coords={},
  		visualization depends on={value \thisrowno{0} \as \PunktI},
  		visualization depends on={value \thisrowno{1} \as \Scheitel},
  		visualization depends on={value \thisrowno{2} \as \PunktII},
  		visualization depends on={value \thisrowno{3} \as \Winkelradius},
  		visualization depends on={value \thisrowno{4} \as \Winkelfarbe},
  		visualization depends on={value \thisrowno{5} \as \Winkelname},
  		visualization depends on={value \thisrowno{6} \as \WinkelExzentrizitaet},
  		nodes near coords style={anchor=center,
  			path picture={
  				\begin{pgfonlayer}{bg}    
  				\draw pic [angle radius=\Winkelradius cm,%
  				fill=\Winkelfarbe!40, draw=\Winkelfarbe,
  				"$\Winkelname$", angle eccentricity =\WinkelExzentrizitaet,
  				text=\Winkelfarbe%
  				] {angle = P\PunktI--P\Scheitel--P\PunktII};
  				\end{pgfonlayer}
  		}},%
  		]
  		table[header=true, x expr=3.37452506686058084640, y expr=0.62932039104983827915]{
  			Punkt1 Scheitel Punkt2 Winkelradius[cm] Winkelfarbe Winkelname WinkelExz
  			2 1 3 0.2 blue!70!black {\alpha} 1.7 \\
  			1 3 4 0.3 blue!70!black {\beta} 1.7 \\
  			4 3 5 0.4 blue!70!black {\gamma} 1.5 \\
  			3 5 6 0.4 blue!70!black {\delta} 1.5 \\
  			4 8 9 0.6 blue!70!black {\epsilon} 1.5 \\
  			6 5 11 0.3 blue!70!black {\zeta} 1.7 \\
  			5 11 12 0.2 blue!70!black {\eta} 1.7 \\
  			9 14 15 1.2 blue!70!black {\theta} 1.2 \\
  			12 11 16 0.4 blue!70!black {\iota} 1.5 \\
  			19 1 2 0.4 blue!70!black {\kappa} 1.5 \\
  			20 19 1 0.4 blue!70!black {\lambda} 1.5 \\
  			25 19 20 0.6 blue!70!black {\mu} 1.4 \\
  			18 16 26 0.4 blue!70!black {\nu} 1.6 \\
  		};
  		\end{axis}
  		
  		\node[below, font=\scriptsize] at (P1) {$P1$};
  		\node[right, font=\scriptsize] at (P2) {$P2$};
  		\node[below, font=\scriptsize] at (P3) {$P3$};
  		\node[above, font=\scriptsize] at (P4) {$P4$};
  		\node[below, font=\scriptsize] at (P5) {$P5$};
  		\node[left, font=\scriptsize] at (P6) {$P6$};
  		\node[left, font=\scriptsize] at (P7) {$P7$};
  		\node[below, font=\scriptsize] at (P8) {$P8$};
  		\node[below, font=\scriptsize] at (P9) {$\quad P9$};
  		\node[above, font=\scriptsize] at (P10) {$P10\quad$};
  		\node[below, font=\scriptsize] at (P11) {$P11$};
  		\node[below, font=\scriptsize] at (P12) {$P12\quad$};
  		\node[right, font=\scriptsize] at (P13) {$P13$};
  		\node[right, font=\scriptsize] at (P14) {$P14$};
  		\node[above, font=\scriptsize] at (P15) {$P15$};
  		\node[left, font=\scriptsize] at (P16) {$P16$};
  		\node[below, font=\scriptsize] at (P17) {$\quad P17$};
  		\node[left, font=\scriptsize] at (P18) {$P18$};
  		\node[right, font=\scriptsize] at (P19) {$P19$};
  		\node[right, font=\scriptsize] at (P20) {$P20$};
  		\node[right, font=\scriptsize] at (P21) {$P21$};
  		\node[left, font=\scriptsize] at (P22) {$P22$};
  		\node[above, font=\scriptsize] at (P23) {$P23$};
  		\node[right, font=\scriptsize] at (P24) {$P24$};
  		\node[right, font=\scriptsize] at (P25) {$P25$};
  		\node[left, font=\scriptsize] at (P26) {$P26$};
  		\node[above, font=\scriptsize] at (P27) {$P27$};
  		\node[above, font=\scriptsize] at (P28) {$P28$};
  		\node[left, font=\scriptsize] at (P29) {$P29$};
  		\node[above, font=\scriptsize] at (P30) {$P30$};
  		\node[below, font=\scriptsize] at (P31) {$P31$};
  		\node[above, font=\scriptsize] at (P32) {$P32$};
  		\node[right, font=\scriptsize] at (P33) {$P33$};
  		\node[right, font=\scriptsize] at (P34) {$P34$};
  		\node[above, font=\scriptsize] at (P35) {$P35$};
  		\node[above, font=\scriptsize] at (P36) {$P36$};
  		\node[below, font=\scriptsize] at (P37) {$\quad P37$};
  		\node[right, font=\scriptsize] at (P38) {$P38$};
  		\node[above, font=\scriptsize] at (P39) {$\quad\quad P39$};	
  		\node[left, font=\scriptsize] at (P40) {$P40$};
  		\node[left, font=\scriptsize] at (P41) {$P41$};
  		\node[below, font=\scriptsize] at (P42) {$P42$};
  		\node[right, font=\scriptsize] at (P43) {$P43$};
  		\node[above, font=\scriptsize] at (P44) {$P44\quad$};
  		\node[right, font=\scriptsize] at (P45) {$P45$};
  		\node[left, font=\scriptsize] at (P46) {$P46$};
  		\node[left, font=\scriptsize] at (P47) {$P47$};
  		\node[right, font=\scriptsize] at (P48) {$P48$};
  		\node[right, font=\scriptsize] at (P49) {$P49$};
  		\node[below, font=\scriptsize] at (P50) {$P50$};
  		\node[left, font=\scriptsize] at (P51) {$P51$};
  		\node[below, font=\scriptsize] at (P52) {$P52$};
  		\node[above, font=\scriptsize] at (P53) {$P53$};
  		\node[below, font=\scriptsize] at (P54) {$P54$};

  		\end{tikzpicture}
  		
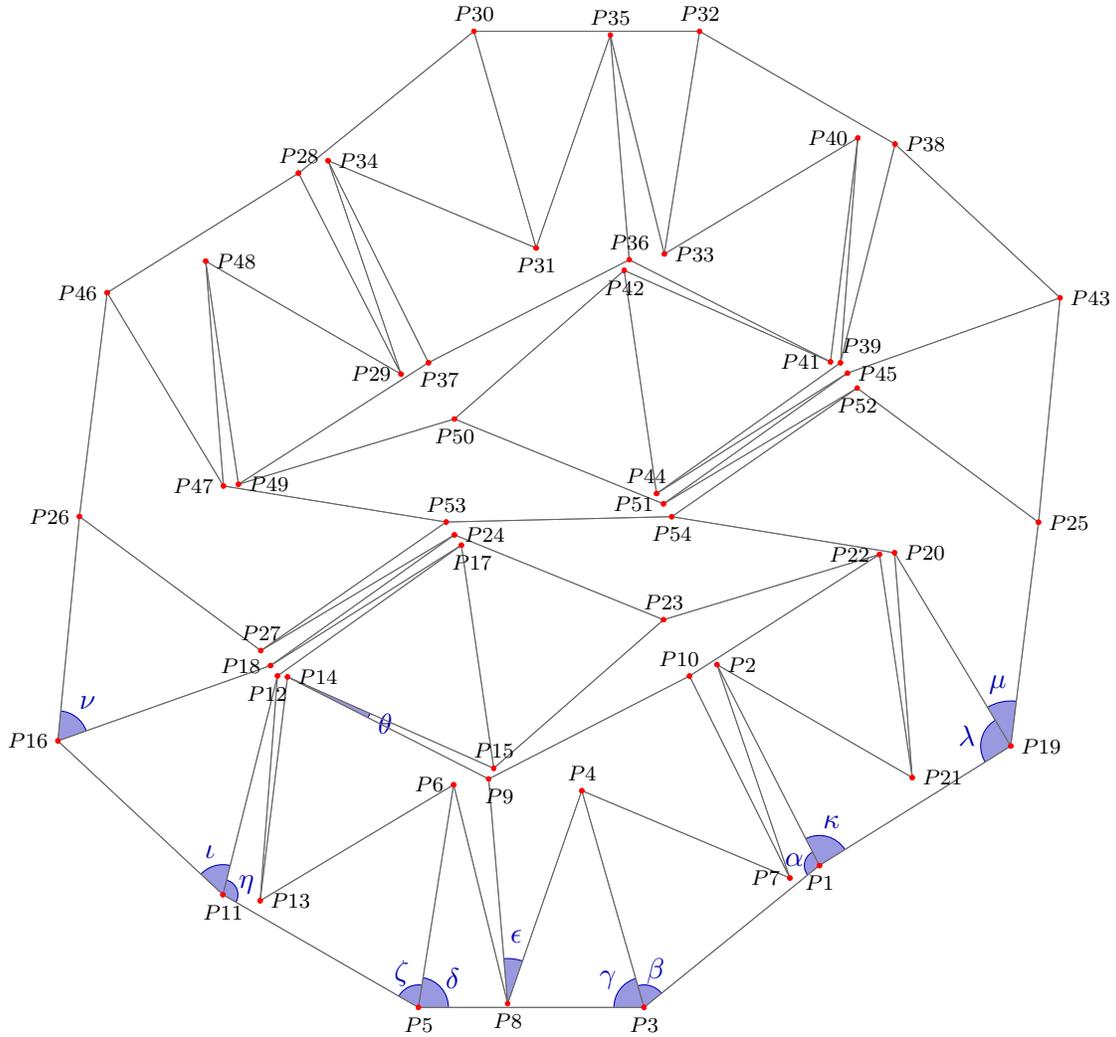
\captionof{figure}{Only known example of a 3-regular matchstick graph of girth 5 consisting of 54 vertices.}
  	\end{minipage}
  \end{center}
  \quad
  
  \begin{center}
  	\Large\textsc{\\Open problems and conjectures about\\matchstick graphs}
  \end{center}
  
  Kurz and Mazzuoccolo write in their publication of 2010: "Our knowledge about matchstick graphs is still very limited. It seems to be hard to obtain rigid mathematical results about them. \cite{Kurz}" Unfortunately, nearly a decade later this statement is still valid. One of the reasons for this, besides the mathematically difficulty, is certainly that still too few people deal with this topic. We hope this article could change this. The following list of open problems and conjectures should be an incentive for the interested reader.
  \\ \\
  \textbf{Problem 1}: Does a 3-regular matchstick graph of girth 5 with less than 54 vertices exists?
  \\ \\
  \textbf{Problem 2}: Examples of 4-regular matchstick graphs are currently known for all number of vertices $\geq$ 52 except for 53, 55, 56, 58, 59, 61, and 62 \cite{Existence}. Try to find an example for one of the missing numbers of vertices, especially with less than 52 vertices.
  \\ \\
  \textbf{Problem 3}:
  Is there an elegant proof or an algorithm which is relatively fast in practice to show that a planar graph with less than 10 vertices can be a matchstick graph?
  \\ \\
  \textbf{Conjecture}: The so-called "Harborth-Graph" with 52 vertices is the smallest possible example of a 4-regular matchstick graph that contains two lines of symmetry, as in a rhombus. A proof of this conjecture may be possible.
  
  \begin{center}
  	\Large\textsc{\\References}
  \end{center}
  
  \begingroup
  \renewcommand{\section}[2]{}
  
  \endgroup

\end{document}